\documentclass{article}
\usepackage{amsmath}
\usepackage{amssymb,latexsym,eucal}
\usepackage{authblk}

\usepackage[english,francais]{babel}

\newtheorem{thm}{Theorem}[section]
 \newtheorem{prop}{Proposition}[section]

\newtheorem{cor}{Corollary}[section]

\newcommand{\ZZ}{\mathbb{Z}}
\newcommand{\NN}{\mathbb{N}}

\def\og{\leavevmode\raise.3ex\hbox{$\scriptscriptstyle\langle\!\langle$~}}
\def\fg{\leavevmode\raise.3ex\hbox{~$\!\scriptscriptstyle\,\rangle\!\rangle$}}

\begin{document}

\title{About a low complexity class of Cellular Automata\footnote{Published in C. R. Acad. Sci. Paris,  Ser I 346 (2008) 995-998 }}

\author[1]{Pierre Tisseur }
\affil[1]{ Centro de Matematica, Computa\c{c}\~ao e Cogni\c{c}\~ao, Universidade Federal do ABC, Santo Andr\'e, S\~ao Paulo, Brasil \thanks{E-mail address: \texttt{pierre.tisseur@ufabc.edu.br}}}
 
\sloppy 
\date{\empty }

\maketitle

\begin{abstract}
Extending to all probability measures the notion of $\mu$-equicontinuous cellular automata  introduced for  Bernoulli measures by Gilman,
 we show that the  entropy is null if $\mu $ is an  invariant measure and that 
the sequence of image measures of  a shift ergodic measure  by iterations of such automata  converges in
Ces\`{a}ro mean to an invariant measure $\mu_c$.
Moreover this cellular automaton  is  still  
$\mu_c$-equicontinuous  and the set of  periodic points is 
dense in the topological support of the  measure $\mu_c$. 
The last property is also true when $\mu$ is invariant and shift ergodic.
 \end{abstract}

\section{Introduction, definitions}

\label{}

Let $A$ be a finite set. We denote by  $A^\ZZ$, the set of bi-infinite sequences
 $x=(x_i)_{i\in\ZZ}$ where $x_i\in A$. 
 We endow $A^\ZZ$ with the product
topology of the discrete topologies on  $A$.  A point $x\in A^\ZZ$ is called a
configuration. The shift $\sigma \colon A^\ZZ\to A^\ZZ$ is defined by :
$\sigma (x)=(x_{i+1})_{i\in \ZZ}$. 
 A cellular automaton  (CA) is a continuous self-map $F$ on
$A^\ZZ$ commuting with the shift. The Curtis-Hedlund-Lyndon theorem 
states  that for every cellular automaton $F$ there exist an integer $r$ and a
block map $f$ from $A^{2r+1}$ to $A$ such that $F(x)_i=f(x_{i-r},\ldots ,x_i
,\ldots ,x_{i+r}).$ 
 The integer $r$ is called the radius of the cellular
automaton. 
For integers $i,j$ with $i\le j$  we denote by $x(i,j)$ the word $x_i\ldots
x_j$ and by $x(i,\infty )$ the infinite sequence $(v_n)_{n\in\NN}$ such that for
all $n\in\NN$ one has $v_n=x_{i+n}$.
For any integer  $n\ge 0$ and point $x\in A^\ZZ$, we denote by  $B_n(x)$ the set of  points $y$ such that for all
$i\in\NN$, one has $F^i(x)(-n,n)=F^i(y)(-n,n)$ and by  $C_n(x)$   
the set of  points $y$ such that $y_j=x_j$ with $-n\le j\le n$.
A point $x\in A^\ZZ$ is called an equicontinuous point if for all positive integer $n$ there exists another
positive integer $m$ such that $B_n(x)\supset C_m(x)$.
A point $x$ is  $\mu$-equicontinuous if  for all $m\in\NN$ one has 
$
\lim_{n\to\infty}\frac{\mu \left(C_n(x)\cap B_m(x)\right )}{\mu (C_n(x))}=1.
$
 In this paper, we  call $\mu$-equicontinuous CA any cellular automaton with  a set of full measure of 
$\mu$-equicontinuous points.
 Clearly an equicontinuous point which belongs to 
 $S(\mu )=\overline{\{x\in A^\ZZ\vert \mu(C_n(x))>0\vert \forall n\in\NN\}}$, (the topological support of 
$\mu$) is also a $\mu$-equicontinuous point. When $\mu$ is 
 a shift ergodic measure, the existence of  $\mu$-equicontinuous points
 implies than the cellular automaton is $\mu$-equicontinuous (see \cite{GI87}).

These definitions was motivated by the work of Wolfram (see \cite{Wo86}) that have proposed a first empirical classification
 based on computer simulations. In \cite{GI87} Gilman introduce a formal and measurable classification by 
  dividing the set of CA in tree parts (CA with equicontinuous points, CA without equicontinuous points but with 
  $\mu$-equicontinuous points, $\mu$-expansive CA).
  The Gilman's classes are defined thanks to a Bernoulli measure not necessarily invariant and corresponds to the 
  Wolfram's simulations based on random entry. 
  Here we study some properties of the $\mu$-equicontinuous class that 
  allows to construct easily invariant measures (see Theorem \ref{thm1}) and we try to describe what kind 
  of dynamic characterizes $\mu$-equicontinuous CA when $\mu$ is an invariant measure. 
  Finally, remark that the comparison between  equicontinuity (see some properties of this class 
in \cite{BLTI} and \cite{TI2000}) and $\mu$-equicontinuity take more sense when we study the restriction of the 
automaton to $S(\mu )$   (see  Section 4 for comments and examples).

\section{Statement of the results}

\subsection{Gilman 's Results}

\begin{prop}\label{pro2}\cite{GI88}
If  $\exists x$ and  $m\neq 0$ such that $B_n(x)\cap \sigma^{-m}B_n(x)\neq\emptyset$
with $n\ge r$ (the radius of the automaton $F$) then the common
sequence { $(F^i(y)(-n,n))_{i\in\NN}$} of all points $y\in B_n(x)$ is ultimately periodic.
\end{prop}
\medskip

In \cite{GI88} Gilman state the following result for any  Bernoulli  measure $\mu$. The proof uses only the shift 
ergodicity of these measures and can be extended to any shift ergodic measure.
\medskip

\begin{prop}\cite{GI88}\label{pro3}
Let $\mu$ be a shift ergodic measure. If a cellular automaton $F$ has a $\mu$-equicontinuous point,
 then for all  $\epsilon >0$ there exists a $F$-invariant closed set $Y$
such that $\mu (Y)>1-\epsilon$ and the restriction of $F$ to $Y$ is equicontinuous. 
\end{prop}

\subsection{New Results}

\begin{prop}\label{pro4}
The measure entropy $h_\mu (F)$ of a $\mu$-equicontinuous and $\mu$-invariant 
cellular automaton $F$ (with $\mu$ not necessarily shift invariant)  is equal to zero.
\end{prop}
\medskip

\begin{prop}\label{pro5}
If a  cellular automaton $F$ has some $\mu$-equicontinuous points where $\mu$ is a $F$-invariant and shift ergodic 
measure then the set of $F$-periodic points is dense in the topological support of $\mu$.
\end{prop}
\medskip

\begin{thm}\label{thm1}
Let $\mu$ be a shift-ergodic measure. If a cellular automaton $F$ has some
$\mu$-equicontinuous points then the sequence 
$(\mu_n)_{n\in\NN}=(\frac{1}{n}\sum_{i=0}^{n-1}\mu\circ F^{-i})_{n\in\NN}$ 
converges vaguely to an invariant measure $\mu_c$.
\end{thm}
\medskip

\begin{thm}\label{thm2}
If $\mu$ is a shift ergodic measure and $F$ a $\mu$-equicontinuous cellular automaton then $F$ is also a 
$\mu_c$-equicon\-ti\-nu\-ous cellular automaton.
\end{thm}
\medskip

\begin{cor}\label{pro8}
If $\mu_c=\lim_{n\to\infty}\frac{1}{n}\sum_{i=0}^{n-1}\mu\circ F^{-i}$ where $\mu$ is a shift ergodic measure 
  and $F$ is a cellular automaton with $\mu$-equicontinuous points   then 
 the set of $F$-periodic points is dense in $S(\mu_c )$. 
\end{cor}

\section{Proofs (Sketches)}

\subsection{Proof of Proposition \ref{pro4}} 

Denote by $(\alpha_p )_{p\in\NN}$ the partition of $A^\ZZ$ by the $2p+1$ central coordinates and 
remark that $h_\mu (F)=\lim_{p\to\infty}h_\mu (F,\alpha_p )$ where  $h_\mu (F,\alpha_p )$ denote 
 the  measurable entropy  with respect to the partition $\alpha_p$.
Using the Shannon-McMillan-Breiman Theorem, we can show that $\forall p\in\NN$, there exists $m\in\NN$ such that 
$h_\mu (F,\alpha_p)\le \int \lim_{n\to\infty}\frac{-\log \mu(B_m(x))}{n}d\mu (x)=0$.

\subsection{Proof of Theorem \ref{thm1}}

It is sufficient to show  
that for all $x\in S(\mu)$ and $m\in\NN$ the sequence $\left(\mu_n (C_m(x))\right)_{n\in\NN}$ converges.
From Proposition \ref{pro3} there exists a set $Y_{\epsilon}$
of measure greater than $1-\epsilon$ such that for all points $y\in Y_{\epsilon}$ and positive integer $k$ 
 the sequences $(F^n(y)(-k,k))_{n\in\NN}$ are eventually periodic with preperiod $pp_\epsilon (k)$ and
period $p_\epsilon (k)$.
We get that  
$\mu_n (C_m(x)\cap Y_{\epsilon})=\frac{1}{n}\sum_{i=0}^{pp_\epsilon (k) -1}\mu \left
( F^{-i} \left
(C_m(x)\right )\cap Y_{\epsilon}\right )$ $
 +\frac{1}{n}\sum_{i=pp_\epsilon (k)}^{n-1}\mu \left ( F^{-i}
\left (C_m(x)\right )\cap Y_{\epsilon}\right )
$ for all $x\in A^\ZZ$ and integer $k\ge m$.
Remark that 
the first term tends to $0$ and the periodicity of the second one implies that  
$
\lim_{n\to\infty}\mu_n (C_m(x)\cap Y_{\epsilon})=\frac{1}{p_\epsilon (k)}
\sum_{i=0}^{p_\epsilon (k)-1}\mu
\left ( F^{-(i+pp_\epsilon (k))} (C_m(x) \cap Y_{\epsilon }\right ).
$
Moreover  we have
$\lim_{\epsilon\to 0}\mu_n(C_m(x)\cap Y_{\epsilon })=\mu_n(C_m(x))$.
Since for all $x$ and $m\in\NN$ one has 
$
\left\vert\mu_n (C_m(x)\cap Y_{\epsilon } )-\mu_n (C_m(x))\right\vert \le \frac{n\epsilon}{n}
=\epsilon
$
 the convergence is uniform with respect to $\epsilon$.
It follows that we can reverse the limits and obtain that  
$$
\mu_c=\lim_{n\to\infty}\frac{1}{n}\sum_{i=0}^{n-1}\mu\circ F^{-i}(C_m(x))
=\lim_{n\to\infty}\frac{1}{n}\sum_{i=0}^{n-1}\lim_{\epsilon\to 0}
\mu\circ F^{-i}(C_m(x)\cap Y_{\epsilon} )
$$
$$
=\lim_{\epsilon\to 0}\lim_{n\to\infty}\frac{1}{n}\sum_{i=0}^{n-1}
\mu\circ F^{-i}(C_m(x)\cap Y_{\epsilon } )
$$
$$
=\lim_{\epsilon\to 0}\frac{1}{p_\epsilon (k)}\sum_{i=0}^{p_\epsilon (k)-1}\mu
\left ( F^{-(i+pp_\epsilon (k))} (C_m(x)
)  \cap Y_{\epsilon }\right )=\mu_c(C_m(x)).
$$
The invariance of converging subsequences of $(\mu_n)_{n\in\NN}$ is a classical result. 
\hfill$\Box$

\subsection{Proof of Proposition \ref{pro5}}
$\mbox{ }$ \vskip -.6 cm
Since $\mu$ is a shift ergodic measure and there exist a $\mu$-equicontinuous points $x$, for all 
 $m\in\NN$ and $z\in S(\mu)$ there exist $(i,j)\in\NN^2$ such that  
$\mu \left(C_p(z)\cap \sigma^{-(i+p)}B_r(x)\cap \sigma^{j+p}B_r(x)=:S\right)>0$ ($r$ is the radius of the CA).
From the Poincar\'{e} recurrence theorem, for all $z\in S(\mu)$, there exists $m\in\NN$ and $y\in S$ such that 
$F^m(y)(-r-p-i,j+p-r-1)=y(-r-p-i,j+p-r-1)$. From the Proof of Proposition \ref{pro2} (see \cite{GI88}), the shift periodic point 
$\overline{w}=\ldots www\ldots$ such that $\overline{w}(-r-p-i,j+p-r-1)=w=y(-r-p-i,j+p-r-1)$ belongs to S and 
 since the $F$ orbit of each $y'\in S\cap\{y''\in A^\ZZ\vert\, y''_l=y_l\vert (-r-p-i\le l\le j+p-r-1) \} $ 
 share the same central coordinates,
it follows that  
$F^m(\overline{w})(-r-p-i,j+p-r-1)=w=\overline{w}(-r-p-i,j+p-r-1)$ which implies that $F^m(\overline{w})=\overline{w}$ and 
permit to conclude.



\subsection{Proof of Theorem \ref{thm2} and Corollary \ref{pro8}}
$\mbox{ }$ \vskip -.6 cm
Let $x$ be a $\mu$-equicontinuous point.  
 For all $m\in\NN$, define $Y_m:=\cup_{i,j\in\NN^2} (\sigma^{-i-m}B_r(x)$ $\cap \sigma^{j+m}B_r(x))$ 
 ($r$ is the radius of $F$) and $\Omega_m=\lim_{n\to\infty}\cap_{j=0}^n\cup_{i=j}^\infty F^i(Y_m)$ 
  (the omega-limit set of $Y_m$ under $F$).
Since $\mu$ is a shift ergodic measure and $\mu (B_r(x))>0$, for all $m\in\NN$,  we get that
$\mu (Y_m)=1$ and consequently $\mu_c(\Omega_m)=1$.
Let $\Lambda (F)$ be the omega-limit set of $A^\ZZ$. 
Using the eventual periodicity of $(F^n(x)(-r,r))_{n\in\NN}$ (see Proposition \ref{pro2}), it can be proved that the 
 omega-limit set of 
 $B_r(x)$ is a finite union of sets 
$B_r(z_l)\cap\Lambda (F) $ ($0\le l\le p-1$).
This implies that 
$\Omega_m=\cup_{z\in [z_0\ldots z_{p-1}]}
\cup_{i,j\in\NN^2} \left(\sigma^{-i-m}B_r(z)\cap\sigma^{j+m}B_r(z)\right)\cap\Lambda(F)$ and 
it follows that  for all $z\in S(\mu_c)$ and $k\in\NN$, the inequality 
$\mu_c\big(C_k(z)\cap\Omega_k\big)>0$ implies that there 
always exist a point $z'$ and integers $i,j\ge m$ such that 
 $\mu_c \big(C_p(z)\cap \sigma^{-(i+p)}$ $ B_r(z')\cap \sigma^{j+p}B_r(z')\big)>0$.
 Using final arguments of the proof of Proposition \ref{pro5}, the last inequality is sufficient to show Corollary \ref{pro8}.
For any measurable set $E$, define 
$E^{\mu_c}=\{y\in E\vert \lim_{n\to\infty}\frac{\mu_c(C_n(y)\cap E)}{\mu_c(C_n(y))}=1\}$.
For all $m\in\NN$, define $\Omega'_m:=\cup_{z\in [z_0\ldots z_{p-1}]}
\cup_{i,j\in\NN^2} \left(\sigma^{-i-m}B_r(z)\cap\sigma^{j+m}B_r(z)\right)^{\mu_c}$ $\cap\Lambda(F)$ and  
 denote by $\Omega$ the set $\cap_{m\in\NN}\Omega'_m$. Since for all measurable set $E$, one has 
 $\mu_c(E^{\mu_c})=\mu_c(E)$,  for all $m\in\NN$, we get that   
$\mu_c(\Omega'_m)=1$ and consequently $\mu_c(\Omega )=1$.
Since  for all $y\in\Omega$ and 
$k\in\NN$ there exist integers $i,j\ge k$ and a point $z'$ such that 
$y\in \sigma^{-i}(B_r(z')\cap \sigma^j B_r(z'))^{\mu_c}$, 
 we obtain that $y\in B_m^{\mu_c}(y)$ which finish the proof.

\section{Example of $\mu$-equicontinuous CA without equicontinuous points}
$\mbox{ }$ \vskip -.6 cm
In \cite{GI87} Gilman gives an example of a $\mu$-equicontinuous CA $F_s$ that has no equicontinuous points.
The automaton $F_s$ act on $\{0,1,2\}^{\ZZ}$ and is defined thank to the following  block map of radius 1:
\small 
 $$
 \begin{array}{|c|c|c|c|c|c|c|c|c|}
 *00&*01&*02&*10&*11&*12&*20&*21&*22\cr
 0&1&0&0&1&0&2&0&2
 \end{array}
 $$
\normalsize 
 
The letter $*$ stands for any letter in $\{0,1,2\}$.
 Considering $0$ as a background element, the 2's move straight down, 1's move to the left 
 and 1 and 2 collide annihilate each other.
In this case the measure $\mu$ is a Bernoulli measure on $\{0,1,2\}^\ZZ$ and the existence of 
$\mu$-equicontinuous points depends on the parameters $p(0),p(1),p(2)$ of this measure.
In \cite{GI87} it is shown that if $p(2)>p(1)$ then the probability that a 2 is never annihilated is positive and 
this implies that there exist $\mu$-equicontinuous points. 
Since the existence or non existence of a sufficient number of 1 in the right side can always modify the central 
coordinates one has $C_m(x)\not\subset B_n(x)$ for all $n,m\in\NN$ which implies that  there is no equicontinuous points. 

Remark that using Theorem \ref{thm1} and \ref{thm2} the automaton $F_s$ is $\mu_c$-equicontinuous if 
$p(2)>p(1)$ but the restriction of $F_s$ to $S(\mu_c)$ always has equicontinuous points 
($S(\mu_c)= \{0,2\}^\NN$ and $F: S(\mu_c)\to S(\mu_c)$ is the identity).
 In \cite{TI2008}, we describe a more complex CA $\mathcal{F}$ such that $\mathcal{F}: S(\mu_c)\to S(\mu_c)$ is $\mu_c$-equicontinous, without
equicontinuous points and the invariant measure $\mu_c$ is construct thanks to Theorem \ref{thm1}. 
$\mbox{ }$ \vskip -5 cm
\section*{Acknowledgements}
\vskip -.3 cm
We whish to acknowledge the support of the CNPq, the Departamento de Matem\'{a}tica da UNESP, S\~{a}o Jos\'{e} 
do Rio Preto e o Centro de Matem\'{atica, Computa\c{c}\~{a}o e Cogni\c{c}\~{a}o}, Universidade Federal do ABC, 
 Brasil in which this work have been done and the referee for his suggestions.
 


\begin{thebibliography}{99}
 
\bibitem{BLTI} F. Blanchard, P. Tisseur, Some properties of cellular automata with equicontinuity
points, Annale de l'intitut Henri Poincar\'e, Probabilit\'es et Statistiques 36(5) (2000) 569-582.
\bibitem{GI87} R. H. Gilman, Classes of linear automata,  Ergodic theory and
Dynamical Systems 7 (1987) 105-118.
\bibitem{GI88} R. H. Gilman, Periodic behaviour of linear automata,  
Dynamical Systems Lecture Note in Mathematics 1342, Springer, New York, (1988) 216-219.
\bibitem{TI2000} P. Tisseur, Cellular automata and Lyapunov exponents, Nonlinearity 13 (2000) 1547-1560.
\bibitem{TI2008} P. Tisseur, Density of periodic points, invariant measures and almost equicontinuous points of 
Cellular automata, Advance in Applied Mathematics 42 (2009) 504-518.
\bibitem{Wo86}
S. Wolfram, Theory and Applications of Cellular Automata, World
Scientific, (1986).

\end{thebibliography}
\end{document}